\documentclass[runningheads]{llncs}
\usepackage[T1]{fontenc}
\usepackage{graphicx}
\usepackage{amsmath}
\usepackage{hyperref}

\begin{document}
%
\title{A database of high precision trivial choreographies for the planar three-body problem}
\titlerunning{High precision trivial choreographies for the planar three-body problem}
%
\author{
I. Hristov\inst{1} \and
R. Hristova\inst{1,2} \and
I. Puzynin\inst{3} \and \\
T. Puzynina\inst{3}\and
Z. Sharipov\inst{3} \and
Z. Tukhliev\inst{3}
}

\authorrunning{I. Hristov, R. Hristova, I. Puzynin et al.}
%

\institute{Faculty of Mathematics and Informatics, Sofia University "St. Kliment Ohridski", Sofia, Bulgaria\\
\email{ivanh@fmi.uni-sofia.bg}\\
\and Institute of Information and Communication Technologies,  Bulgarian Academy of Sciences, Sofia, Bulgaria
\and Meshcheryakov Laboratory of Information Technologies, Joint Institute for Nuclear Research, Dubna, Russia\\
\email{zarif@jinr.ru}}

\maketitle              
\begin{abstract}
  Trivial choreographies are special periodic solutions of the planar three-body problem.
  In this work  we use a modified Newton's method based on the continuous analog of  Newton's method and a high precision arithmetic
  for a specialized  numerical search for new trivial choreographies.
  As a result of the search we computed a high precision database of
  462 such orbits, including 397 new ones. The initial conditions and the periods of  all found solutions are given with 180 correct decimal digits.
  108 of the choreographies are linearly stable, including 99 new ones.
  The linear stability is tested by a high precision computing of the eigenvalues of the monodromy matrices.

\end{abstract}
\section{Introduction}

 A choreography is a periodic orbit in which the three bodies move along one and the same trajectory with a time delay of $T/3$,
 where $T$ is the period of the solution. A choreography is called trivial if it is a satellite (a topological power) of the famous figure-eight choreography
 \cite{Moore,Montgomery}. Trivial choreographies are of special interest because many of them are expected to be stable like the figure-eight orbit.
 About 20 new trivial choreographies with zero angular momentum and bodies with equal masses  are found in \cite{Shuv1,Shuv2,Dimi}, including one new linearly stable choreography,
 which was the first found linearly stable choreography after the famous figure-eight orbit. Many three-body choreographies (345) are also found in \cite{Simo}, but they are with nonzero angular momentum and an undetermined topological type.

 In our recent work \cite{Hundreds}  we made a purposeful (on a small domain of initial conditions) numerical search for figure-eight satellites (not necessarily choreographies).
 For numerical search we used a modification of Newton's method with a larger domain of convergence. The three-body problem is well known with the sensitive dependence on the
 initial conditions. To overcome the obstacle of dealing with this sensitivity and to follow the trajectories correctly for a long time, we used as an ODE solver the high order multiple precision Taylor series method
 \cite{Barrio,Liao1,Liao2}.
 As a result we found over 700 new satellites with periods up to 300 time units, including 45 new choreographies.
 7 of the newly found choreographies are shown to be linearly stable, bringing the number of the known linearly stable choreographies up to 9.

 This work can be regarded as a continuation of our recent work \cite{Hundreds}.
 Now we make a specialized numerical search for new trivial choreographies by using the permuted return proximity condition proposed in \cite{Shuv2}.
 We consider the same searching domain and the same searching grid step as in  \cite{Hundreds}. Considering pretty long periods (up to 900 time units), which are much longer
 than those in the previous research, allows us to compute a high precision database of 462 trivial choreographies, including 397 new ones. 99 of the newly found choreographies are linearly stable,
 so the number of the known linearly stable choreographies now rises to 108.

\section{Differential equations describing the bodies motion}
The bodies are with equal masses and they are treated as point masses.
A planar motion of the three bodies is considered.
The normalized differential equations describing the motion of the bodies are:
\begin{equation}
\ddot{r}_i=\!\!\!\sum_{j=1,j\neq i}^{3} \frac{(r_j-r_i)}{{\|r_i -r_j\|}^3}, \quad i=1,2,3.
\end{equation}
The vectors $r_i$, $\dot{r}_i$ have two components: $r_i=(x_i, y_i)$, $\dot{r}_i=(\dot{x}_i, \dot{y}_i)$.
The system (1) can be written as a first order one this way:
\begin{equation}
\dot{x}_i={vx}_i, \hspace{0.15 cm} \dot{y}_i={vy}_i, \hspace{0.15 cm} \dot{vx}_i=\!\!\!\sum_{j=1,j\neq i}^{3} \frac{(x_j-x_i)}{{\|r_i -r_j\|}^3}, \hspace{0.15 cm} \dot{vy}_i=\!\!\!\sum_{j=1,j\neq i}^{3} \frac{(y_j-y_i)}{{\|r_i -r_j\|}^3}, \quad i=1,2,3
\end{equation}
We solve numerically the problem in this first order form. Hence we have a vector of 12 unknown functions
$ X(t)={(x_1, y_1, x_2, y_2, x_3, y_3, {vx}_1, {vy}_1, {vx}_2, {vy}_2, {vx}_3, {vy}_3)}^\top$.
Let us mention that this first order system actually coincides with the Hamiltonian formulation of the problem.

We search for periodic planar collisionless orbits as in \cite{Shuv1,Shuv2}: with zero angular momentum and symmetric initial
configuration with parallel velocities:
\begin{equation}
\begin{aligned}
(x_1(0),y_1(0))=(-1,0), \hspace{0.2 cm} (x_2(0),y_2(0))=(1,0), \hspace{0.2 cm} (x_3(0),y_3(0))=(0,0) \\
({vx}_1(0),{vy}_1(0))=({vx}_2(0),{vy}_2(0))=(v_x,v_y) \hspace{2 cm}\\
({vx}_3(0),{vy}_3(0))=-2({vx}_1(0),{vy}_1(0))=(-2v_x, -2v_y) \hspace{1.5 cm}
\end{aligned}
\end{equation}
The velocities $v_x\in [0,1], v_y\in [0,1]$ are parameters. We denote the periods of the orbits with $T$. So, our goal is to find triplets $(v_x, v_y, T)$
for which the periodicity condition $X(T)=X(0)$ is fulfilled.

\section{Numerical searching procedure}

The numerical searching procedure consists of three stages.
During the first stage we search for candidates for correction with the modified Newton's method, i.e. we compute initial approximations of the triplets $(v_x, v_y, T)$.
In what follows we will use the same notation for $v_x, v_y, T$ and their approximations.
We introduce a square 2D searching grid with stepsize 1/4096 for the parameters $v_x, v_y$ in the same searching domain as in \cite{Hundreds}
(the domain will be shown later in Section 5).
We simulate the system (2) at each grid point $(v_x, v_y)$  up to a prefixed time $T_0=300$.
For an ODE solver we use  the high order multiple precision Taylor series method \cite{Barrio,Liao1,Liao2} with a variable stepsize strategy from \cite{Jorba}.
Because we concentrate only on choreographies (that is why we call this search specialized),
we take as in \cite{Shuv2} the candidates to be the triplets $(v_x, v_y, T)$ for which the cyclic permutation return proximity function $R_{cp}$ (the minimum):
$$R_{cp}(v_x,v_y,T_0) = \min_{1<t\leq T_0} {\|\hat{P}X(t)-X(0)\|}_{2}$$
is obtained at $t=T/3$ and is less than 0.1.
We also set the constraint that $R_{cp}$ has a local minimum on the grid for $v_x, v_y$.
Here $\hat{P}$ is a cyclic permutation of the bodies' indices. Using cyclic permutation return proximity instead of the standard return proximity,
reduces three times the needed integration time at the first stage (now we obtain candidates with periods $T$  up to 900).

During the second stage we apply the modified Newton's method, which has a larger domain of convergence
than the classic Newton's method. Convergence during this stage means that a periodic orbit is found.
The following linear algebraic system with a $12\times3$ matrix for the corrections $\Delta v_x,
\Delta v_y, \Delta T$ has to be solved at each iteration step \cite{Abad}.
\begin{equation}
\begin{pmatrix}
\frac{\partial x_1}{\partial v_x}(T) & \frac{\partial x_1}{\partial v_y}(T) & \dot{x}_1(T)\\
\frac{\partial y_1}{\partial v_x}(T) & \frac{\partial y_1}{\partial v_y}(T) & \dot{y}_1(T)\\
\frac{\partial x_2}{\partial v_x}(T) & \frac{\partial x_2}{\partial v_y}(T) & \dot{x}_2(T)\\
\frac{\partial y_2}{\partial v_x}(T) & \frac{\partial y_2}{\partial v_y}(T) & \dot{y}_2(T)\\
\frac{\partial x_3}{\partial v_x}(T) & \frac{\partial x_3}{\partial v_y}(T) & \dot{x}_3(T)\\
\frac{\partial y_3}{\partial v_x}(T) & \frac{\partial y_3}{\partial v_y}(T) & \dot{y}_3(T)\\
\frac{\partial {vx}_1}{\partial v_x}(T) -1 & \frac{\partial {vx}_1}{\partial v_y}(T) & \dot{vx}_1(T)\\
\frac{\partial {vy}_1}{\partial v_x}(T) & \frac{\partial {vy}_1}{\partial v_y}(T) - 1 & \dot{vy}_1(T)\\
\frac{\partial {vx}_2}{\partial v_x}(T) -1 & \frac{\partial {vx}_2}{\partial v_y}(T) & \dot{vx}_2(T)\\
\frac{\partial {vy}_2}{\partial v_x}(T) & \frac{\partial {vy}_2}{\partial v_y}(T) - 1 & \dot{vy}_2(T)\\
\frac{\partial {vx}_3}{\partial v_x}(T) +2 & \frac{\partial {vx}_3}{\partial v_y}(T) & \dot{vx}_3(T)\\
\frac{\partial {vy}_3}{\partial v_x}(T) & \frac{\partial {vy}_3}{\partial v_y}(T) +2 & \dot{vy}_3(T)
\end{pmatrix}
\begin{pmatrix}
\Delta v_x \\
\Delta v_y \\
\Delta T
\end{pmatrix}
=
\begin{pmatrix}
x_1(0) - x_1(T) \\
y_1(0) - y_1(T) \\
x_2(0) - x_2(T) \\
y_2(0) - y_2(T) \\
x_3(0) - x_3(T) \\
y_3(0) - y_3(T) \\
{vx}_1(0) - {vx}_1(T)\\
{vy}_1(0) - {vy}_1(T)\\
{vx}_2(0) - {vx}_2(T)\\
{vy}_2(0) - {vy}_2(T)\\
{vx}_3(0) - {vx}_3(T)\\
{vy}_3(0) - {vy}_3(T)
\end{pmatrix}
\end{equation}
For the classic Newton's method we correct to obtain the next approximation this way:
$$v_x := v_x + \Delta v_x,\hspace{0.2 cm} v_y := v_y + \Delta v_y, \hspace{0.2 cm} T := T + \Delta T$$
For the modified Newton's method based on the continuous analog of Newton's method \cite{CANM}
we introduce a parameter $\tau_k: 0 < \tau_k \leq 1$, where $k$ is the number of the iteration.
Now we correct this way:
$$v_x := v_x + \tau_k\Delta v_x, \hspace{0.2 cm} v_y := v_y + \tau_k\Delta v_y, \hspace{0.2 cm} T := T + \tau_k\Delta T$$
Let $R_{k}$ be the residual ${\|X(T)-X(0)\|}_2$ at the $k$-th iteration. With a given $\tau_0$, the next $\tau_k$ is computed
by the following adaptive algorithm \cite{CANM}:
\begin{equation}
\tau_k = \left \{
               \begin{array}{ll}
                \min(1,\ \tau_{k-1} R_{k-1} / R_k), &
                 R_k \leq R_{k-1}, \\\\
                \max(\tau_0,\ \tau_{k-1} R_{k-1} / R_k), &
                 R_k > R_{k-1},
               \end{array}
                                  \right .
\end{equation}
The value $\tau_0=0.2$ is chosen in this work. We iterate until the value $R_k$ at some iteration becomes less than some tolerance
or the number of the iterations becomes greater than some number $maxiter$ to detect
divergence. The modified Newton's method has a larger domain of convergence than the classic
Newton's method, allowing us to find more choreographies for a given search grid. To compute the matrix elements in (4),
a system of 36 ODEs (the original 12 differential equations plus the 24 differential equations for the partial derivatives with respect
to the parameters $v_x$ and $v_y$) has to be solved. The high order multiple precision Taylor series method is used again for solving this
36 ODEs system. The linear algebraic system (4) is solved in linear least square sense using QR decomposition based on Householder reflections \cite{Demmel}.

During the third stage we apply the classic Newton's method with a higher precision in order to specify the solutions with more correct digits (180 correct
digits in this work). This stage can be regarded as some verification of the found periodic solutions since we compute the initial conditions and the periods with
many correct digits and we check the theoretical quadratic convergence of the Newton's method.

During the first and the second stage the high order multiple precision Taylor series method with order 154 and 134 decimal digits of precision is used.
During the third stage we make two computations. The first computation is with 242-nd order Taylor series method and 212
digits of precision and the second computation is for verification - with 286-th order method and 250 digits of precision.
The most technical part in using the Taylor series method is the computations of the derivatives for the Taylor's formula, which are done
by applying the rules of automatic differentiation \cite{Barrio2}. We gave all the details for the Taylor series method,
particularly we gave all needed formulas, based on the rules of the automatic differentiation in our work \cite{Ourpaper}.

\section{Linear stability investigation}
The linear stability of a given periodic orbit $X(t), 0 \leq t \leq T,$ is determined by the eigenvalues $\lambda$ of the $12 \times 12 $ monodromy matrix $M(X;T)$ \cite{Roberts}:
$$M_{ij}[X;T] = \frac{\partial X_i(T)}{\partial X_j(0)}, \quad  M(0)=I $$
The elements $M_{ij}$ of the monodromy matrix $M$ are computed in the same way as the partial derivatives in the system (4)
- with the multiple precision Taylor series method using the rules of automatic differentiation (see for details \cite{Ourpaper}).
The eigenvalues of $M$ come in pairs or quadruplets: $(\lambda, \lambda^{-1}, \lambda^{*}, \lambda^{*-1})$. They are of four types:

1) Elliptically stable - $\lambda = \exp(\pm 2 \pi i \nu),$ where $\nu > 0$ (real) is the stability angle. In this case the eigenvalues are on the unit circle.
Angle $\nu$ describes the stable revolution of adjacent trajectories around a periodic orbit.

2) Marginally stable - $\lambda = \pm 1.$

3) Hyperbolic - $\lambda = \pm\exp(\pm \mu),$ where $\mu>0$ (real) is the Lyapunov exponent.

4) Loxodromic - $\lambda = \exp(\pm\mu \pm i \nu),$  $\mu, \nu$ (real)

Eight of the eigenvalues of $M$ are equal to $1$ \cite{Roberts}. The other four determine the linear stability.
Here we are interested in elliptically stable orbits, i.e. the four eigenvalues to be $\lambda_j= \exp(\pm 2 \pi i \nu_j), \nu_j>0, j=1,2$.
For computing the eigenvalues we use a Multiprecision Computing Toolbox \cite{Advanpix} for MATLAB \cite{Matlab}.
First the elements of $M$ are obtained with 130 correct digits and then two computations  with 80 and 130 digits of precision are made
with the toolbox. The four eigenvalues under consideration  are verified by a check for matching the first 30 digits of them and the corresponding
condition numbers obtained by the two computations (with 80 and 130  digits of precision).

\section{Numerical results}

To classify the periodic orbits into topological families we use a topological method from \cite{Mont}.
Each family corresponds to a different conjugacy  class of the free group on two letters $(a,b)$.
Satellites of figure-eight correspond to free group elements ${(abAB)}^k$
for some natural power $k$. For choreographies the power $k$ can not be divisible by 3 (see \cite{Shuv1}).
We use "the free group word reading algorithm" from \cite{Shuv3} to obtain the free group elements.
Together with the triplet $(v_x, v_y, T)$, we compute the scale-invariant period $T^*$.
$T^*$ is defined as $T^{*}=T{|E|}^{\frac{3}{2}},$ where  $E$ is the energy of our initial configuration: $E=-2.5 + 3({v_x}^2+{v_y}^2)$.
Equal $T^*$ for two different initial conditions means two different representations of the same solution (the same choreography).
Some of the choreographies are presented by two different initial conditions.

We found 462 trivial choreographies in total (including Moore's figure-eight orbit and old choreographies, and counting different initial conditions as different solutions).
397 of the solutions are new (not included in \cite{Shuv1,Shuv2,Dimi,Hundreds}).
For each found solution we computed the power $k$ and the four numbers $(v_x,v_y,T,T^*)$ with 180 correct digits.
This data can be seen in \cite{rada3body} together with the plots of the trajectories in the real $x-y$ plane.

As a result of computing the eigenvalues of the monodromy matrices we obtain that 108 of the choreographies are linearly stable (99 new ones).
These 108 choreographies correspond to 150 initial conditions (some of the choreographies are presented by two different initial conditions and have the same $T^*$).
All stability angles $\nu_{1,2}$ with 30 correct digits for the linearly stable solutions
are given in a table in \cite{rada3body}. The distribution of the initial condition points
can be seen in Fig.~1 (the black points are the linearly stable solutions, the white points -- the unstable ones (more precisely, not confirmed to be linearly stable)).
The searching domain is the same as in \cite{Hundreds} and consists of the rectangle $[0.1, 0.33]\times [0.49, 0.545]$ plus the domain with the curved boundary. The three points out of the searching domain
are previously found solutions.

\begin{figure}
\begin{center}
\includegraphics[scale=0.57]{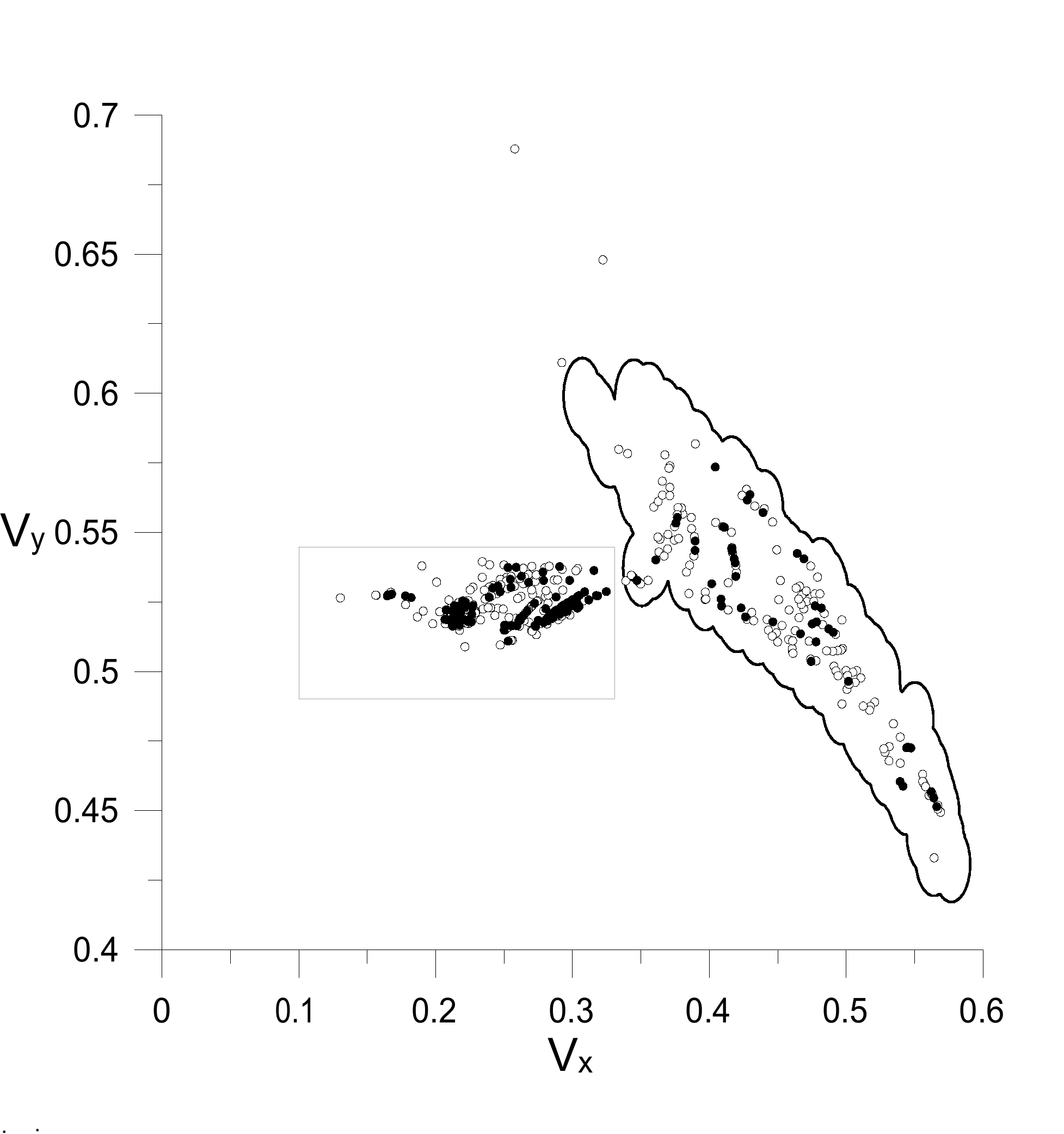}
\caption{Initial velocities for all found solutions, black points are the linearly stable solutions,
white points - unstable (more precisely, not confirmed to be linearly stable)}
\label{fig1}
\end{center}
\end{figure}

Analyzing the linear stability data, we observe 13 pairs of solutions with the same power $k$ and very close $T^*$, where one is linearly stable and the other is of hyperbolic-elliptic type. The pairs have the following property:
"The stability angle of the elliptic eigenvalues of the hyperbolic-elliptic type solution is very close to the larger stability angle of the linearly stable solution, the smaller stability angle of the linearly stable solution is close to zero and the larger hyperbolic eigenvalue of the hyperbolic-elliptic type solution is greater than one but very close to one".
For example, there exists a pair of solutions with {\footnotesize $\nu_1=0.255011944221133753875666925693$}, {\footnotesize $\nu_2=2.19223274459622941216216635818e-05$} and {\footnotesize $\nu=0.255011941995861150357102898351$}, {\footnotesize $\lambda=1.00013775153254718967585223182$}. The initial conditions and the periods for this pair can
be seen in Table 1. The solutions are those with numbers 119 and 120 in \cite{rada3body}.  The trajectories of the three bodies in the real $x-y$ plane can be seen in Fig.~2 and Fig.~3.
A table with the linear stability data for all the 13 pairs can be seen in \cite{rada3body}.

\begin{table}
 \begin{tabular}{ p{0.5cm} p{2.75cm} p{2.75cm} p{2.75cm} p{2.75cm} p{0.3cm} }
 \hline
 $N$ & $\hspace{1.1cm}v_x$ & $\hspace{1.1cm}v_y$ & $\hspace{1.2cm}T$ & $\hspace{1.2cm}T^*$ & $\hspace{0.05cm}k$ \\
 \hline
\tiny{119} & \tiny{0.41817368353651279e0} & \tiny{0.54057212735770067e0} & \tiny{0.52133539095545824e3} & \tiny{0.600424230253006803e3} & \tiny{65} \\
\tiny{120} & \tiny{0.26562094559259036e0} & \tiny{0.5209803403964781e0} & \tiny{0.33548942966568876e3}   & \tiny{0.600424230253006829e3} & \tiny{65} \\
\end{tabular}

{\caption {Data with 17 correct digits for a pair of linearly stable and hyperbolic-elliptic solutions}}
\end{table}

\vspace{-0.5 cm}

\begin{figure}
\begin{center}
\includegraphics[scale=0.49]{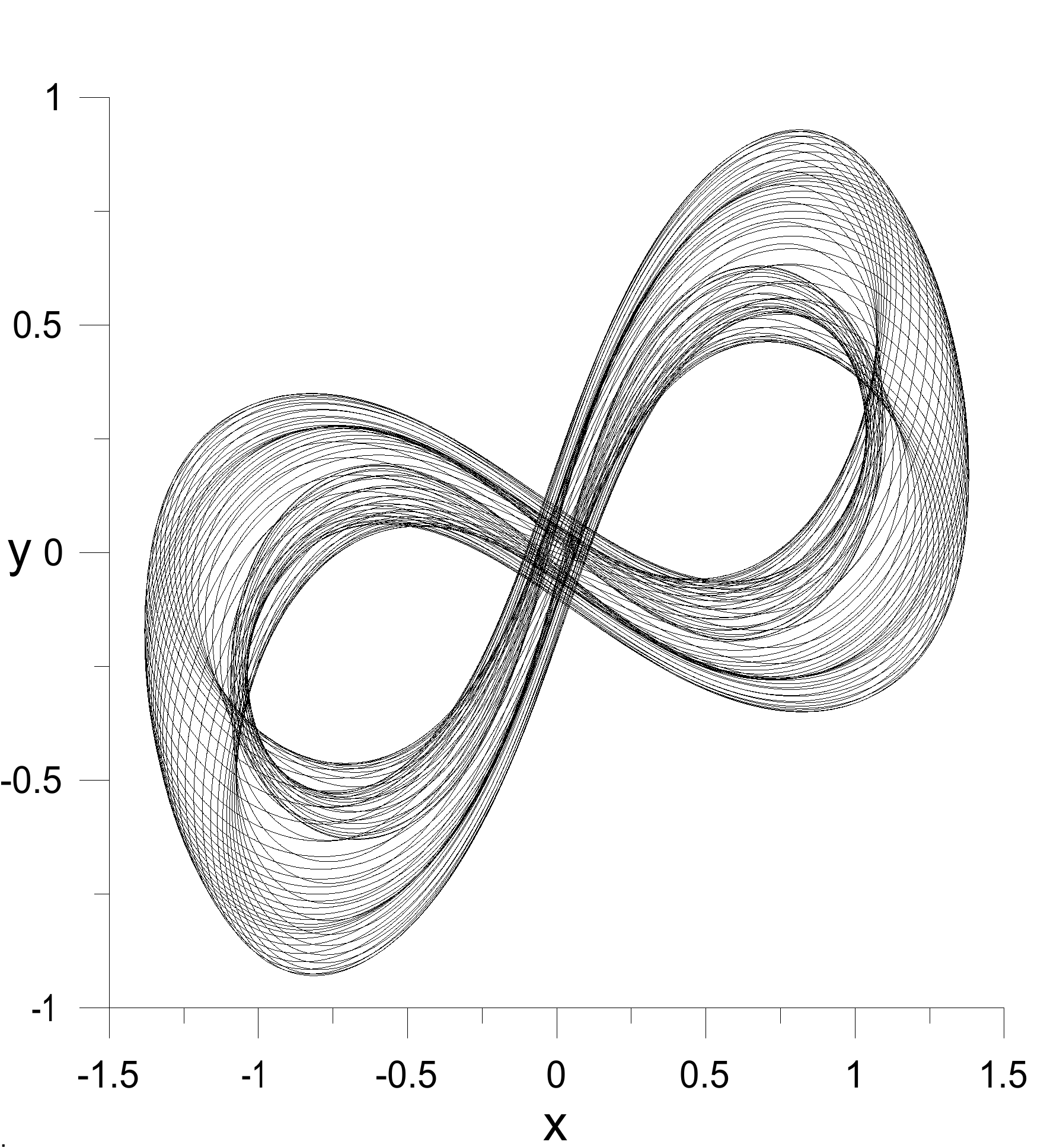}
\caption{Linearly stable solution (solution 119 from Table 1)}
\label{fig:1}
\end{center}
\end{figure}

\begin{figure}
\begin{center}
\includegraphics[scale=0.49]{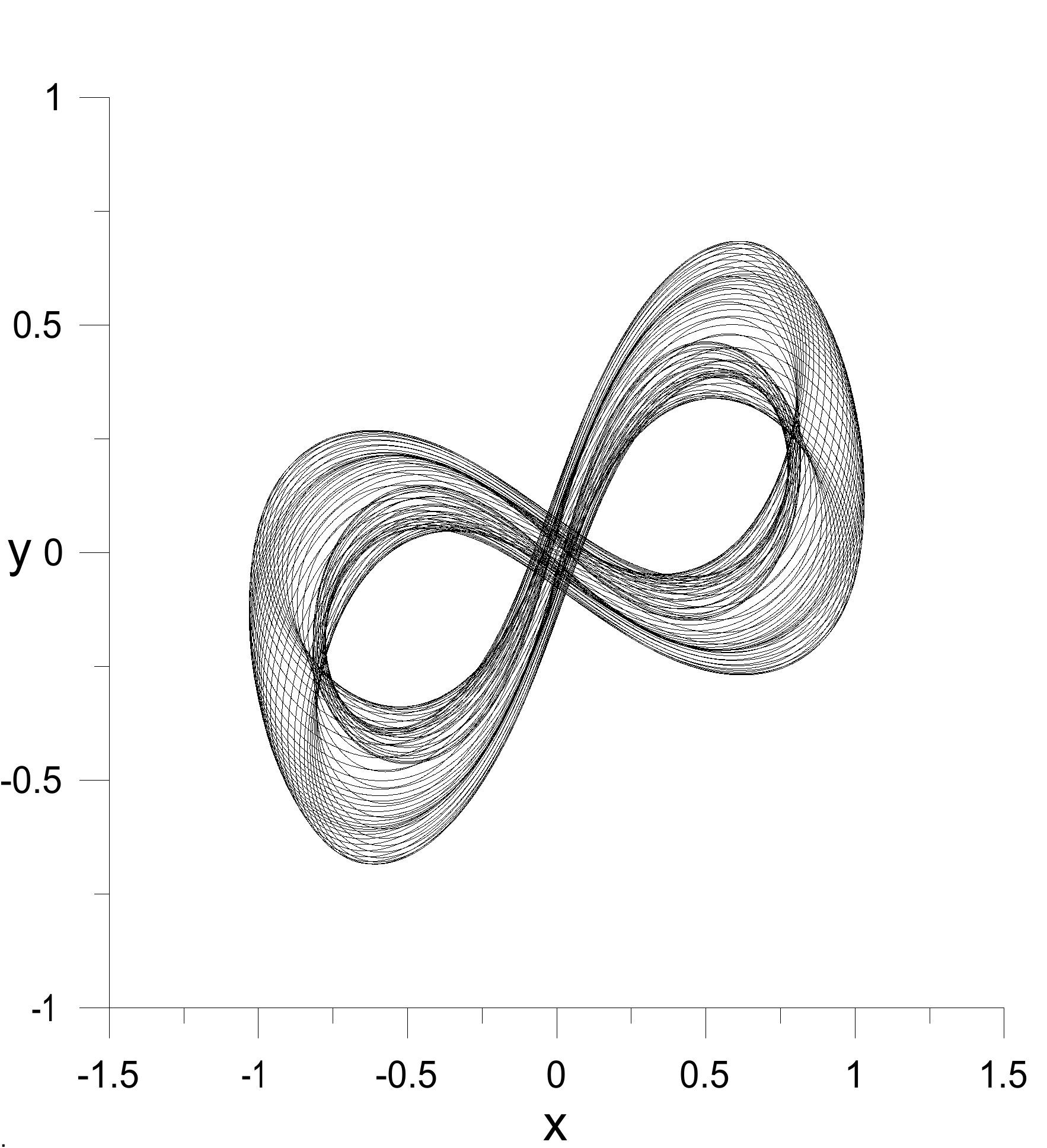}
\caption{Hyperbolic-elliptic type solution (solution 120 from Table 1)}
\label{fig:1}
\end{center}
\end{figure}

The extensive computations for the numerical search are performed in "Nestum" cluster, Sofia, Bulgaria \cite{nestum},
where the GMP library (GNU multiple precision library) \cite{GMP} for multiple precision floating point arithmetic is installed.

\section{Conclusions}
A  modified Newton's method with high precision is successfully used for a specialized numerical search
of new trivial choreographies for the planar three-body problem. Considering pretty long periods (up to 900) allows us to compute a high precision database of 462 solutions (397 new ones).
99 of the newly found solutions are linearly stable, bringing the number of the known linearly stable trivial choreographies from 9 to 108.

\subsubsection{Acknowledgements} We greatly thank for the opportunity to use the computational resources of the "Nestum" cluster, Sofia, Bulgaria.
We would also like to thank Veljko Dmitrashinovich from Institute of Physics, Belgrade University, Serbia for a valuable e-mail discussion and advice,
and his encouragement to continue our numerical search for new periodic orbits.

%
%
%
%

\end{document}